\documentclass[12pt]{article}
\usepackage{amsthm,amsfonts,amssymb,amsmath}
\usepackage{fullpage}
\usepackage{enumerate}
\usepackage[colorlinks=true]{hyperref}
\newcommand{\CL}{{\mathcal {L}}}
\newcommand{\CM}{{\mathcal {M}}}

\newcommand{\CO}{{\mathcal {O}}}

\newcommand{\CX}{{\mathcal {X}}}

\newcommand{\id}{{\mathrm{id}}}

\newcommand{\End}{{\mathrm{End}}}

\newcommand{\NS}{{\mathrm{NS}}}

\newcommand{\Pic}{\mathrm{Pic}}

\DeclareMathOperator{\Spec}{Spec}

\newcommand{\ol }{\overline}

\newcommand{\lra}{\longrightarrow}

\newcommand{\CC}{\mathbb{C}}
\newcommand{\RR}{\mathbb{R}}
\newcommand{\ZZ}{\mathbb{Z}}
\newcommand{\QQ}{\mathbb{Q}}

\newtheorem{thm}{Theorem}[section]

\newtheorem{theorem}[thm]{Theorem}
\newtheorem{proposition}[thm]{Proposition}
\newtheorem{lemma}[thm]{Lemma}

\theoremstyle{definition}

\theoremstyle{remark}

\newcommand\Nef{\mathrm{Nef}}
\newcommand\Amp{\mathrm{Amp}}
\newcommand\Psef{\mathrm{Psef}}
\newcommand\BBig{\mathrm{Big}}
\newcommand\pef{pseudo-effective\ }

\begin{document}

\title{Nef cone and successive minima: an example}
\author{Ruoyi Guo, Xinyi Yuan}
\maketitle

\setcounter{tocdepth}{4}
\tableofcontents

\section{Introduction}

In this paper, we compute the nef cone and the pseudo-effective cone of $C\times J$ for a smooth projective curve $C$ and its Jacobian variety $J$ such that $C\times J$ has the minimal Picard number.
As a consequence, we also compute the successive minima of a height function for the relative setting $C\times J\to J$, and our result shows that Zhang's theorem of  successive minima does not hold for higher dimensional function fields.

\subsection{The cones}

By a \emph{variety} over a field, we mean an integral scheme, separated and of finite type over the field. By a \emph{curve}, we mean a variety of dimension 1. 

Let us review the classical notions of nef cones and the \pef cones. We refer to Lazarsfeld \cite{Laz04} for more details. Let $X$ be a smooth projective variety over a field $k$. Let 
$\mathrm{NS}(X)=\Pic(X)/\Pic^0(X)$
be its N\'eron-Severi group. It is well-known that $\NS(X)$ is a finitely generated abelian group, whose rank is called the \emph{Picard number} $\rho(X)$ of $X$.

Denote $\mathrm{NS}(X)_\RR=\NS(X) \otimes_\ZZ\RR$. 
Denote by 
$$\Amp(X),\quad \Nef(X), \quad \BBig(X), \quad \Psef(X)$$
the subsets of $\xi\in \NS(X)_\RR$ which are ample, nef, big, and pseudo-effective on $X$ respectively. 
They are respectively called the \emph{ample cone}, the \emph{nef cone}, the \emph{big cone}, and the \emph{pseudo-effective cone} of $X$.  
All of them are convex cones in $\NS(X)_\RR$. 
Moreover, $\Nef(X)$ is the closure of $\Amp(X)$, and $\Amp(X)$ is the interior of $\Nef(X)$;  $\Psef(X)$ is the closure of $\BBig(X)$, and $\BBig(X)$ is the interior of $\Psef(X)$. 

It is known that an element $\xi\in \NS(X)_\RR$ is \pef if and only if $\xi + \BBig(X)\subseteq \BBig(X)$. In particular, if a line bundle $L$ on $X$ is effective in that $\Gamma(X,L)\neq 0$, then it is also pseudo-effective. 
In practice, it is usually very difficult to determine these cones for a general variety.
In the following, we will provide an explicit result for a very special example.  

Let $C$ be a smooth projective curve  of genus $g>1$ over an algebraically closed field $k$. Let $J$ be the Jacobian variety of $C$. 
Then the Picard number $\rho(C\times J)\geq 3$ by the classical decomposition 
$$
\NS(C\times J) = \NS(C)\oplus \NS(J) \oplus \End(J).
$$
We refer to 
Theorem \ref{decomposition} for a proof of this decomposition.

It is also easy to construct 3 linearly independent elements of $\NS(C\times J)$.
We will introduce them under our notation explicitly, which will be used throughout this paper. 
Take a line bundle $\alpha\in \Pic(C)$ such that $(2g-2)\alpha=\omega_{C/k}$  in $\Pic(C)$. Denote the Abel--Jacobi map
\[
i_\alpha:C\longrightarrow J,\quad x\longmapsto x-\alpha.
\]
Let $\theta$ be the image of the morphism
\[
C^{g-1}\longrightarrow J,\quad (x_1,\cdots,x_{g-1})\longmapsto x_1+\cdots+x_{g-1}-(g-1)\alpha.
\]
It is well-known that $\theta$ is an ample divisor on $J$, which gives a principal polarization 
\[
\phi:J\stackrel{\thicksim}\longrightarrow J^\vee,\quad  x\longmapsto T^*_{x}\theta-\theta.
\]
Here $T_x:J\to J$ denotes the morphism given by translation by $x$ for $x\in J(k)$.
Let $P$ be the pull-back of the usual Poincar\'e bundle on $J\times J^\vee$ via 
\[
\id\times \phi:J\times J\longrightarrow J\times J^\vee.
\]
Then we have a Poincar\'e line bundle
\[
Q=(i_\alpha\times \id)^*P
\]
 on $C\times J$.

Denote by $p_1:C\times J\to C$ and $p_2:C\times J\to J$ the projections.  
Then $\NS(C\times J)$ contains 3 linearly independent elements
$$
 \alpha_1:= p_1^*\alpha, \quad \theta_2:=p_2^*\theta,\quad Q.
$$
Each of them lies in a component in the above decomposition of  
$\NS(C\times J)$.
In fact, $Q$ corresponds to $1\in \End(J)$, as we can see from Theorem \ref{decomposition}. 
 
In our main results, we will assume the minimal case $\rho(C\times J)=3$. 
We remark that there are plenty of curves $C$ with $\rho(C\times J)=3$. 
In fact, if $C$ is a very general smooth projective curve over $\CC$ in the sense that it is represented by a complex point of the coarse moduli scheme $\CM_g$ over $\CC$ outside a countable union of proper algebraic subsets, then  
$\rho(J)=1$ by a classical theorem of Lefschetz
(cf. \cite{Pir88}) and $\End(J)=\ZZ$ by a classical theorem of Hurwitz
(cf. \cite[\S3]{Cil89} and \cite{Zar00}). 
For such complex curves $C$, we have $\rho(C\times J)=3$. 

The following is our first theorem, which describes the cones of $C\times J$ explicitly.  
\begin{theorem}[Theorem \ref{cone thm}]
 \label{first main thm}
Assume that $\rho(C\times J)=3$. Then
	\[
	\Amp(C\times J)=\BBig(C\times J)=\{a\cdot \alpha_1+b\cdot \theta_2+c\cdot Q:a>0,\ b> 0, \ ab> gc^2\}
	\]
and
	\[
	\Nef(C\times J)=\Psef(C\times J)=\{a\cdot \alpha_1+b\cdot \theta_2+c\cdot Q:a\geq0,\ b\geq 0, \ ab\geq gc^2\}.
	\]
\end{theorem}

The main idea of our proof of the theorem is to construct all nef line bundles on $C\times J$ by geometric operations. In fact, for $m,n\in \ZZ$, consider the morphism 
$$
f_{m,n}: C\times J\longrightarrow J,\quad (x,y)\longmapsto m(x-\alpha)+ny. 
$$
As proved in Proposition \ref{pullback}, we have the pull-back
\[
f_{m,n}^*\theta=gm^2\cdot \alpha_1+n^2\cdot \theta_2+mn\cdot Q.
\]
It is nef on $C\times J$, since $\theta$ is ample on $J$.
This essentially gives all the line bundles on the boundary of the nef cone.

\subsection{Successive minima}

To introduce our second result, let us first review a geometric version of Zhang's theorem of successive minima. For simplicity, we restrict to curves over function fields.

Let $B$ be a projective variety  of dimension $d>0$ over a field $k$, and let $K=k(B)$ be the function field of $B$. Let $X$ be a projective curve over $K$. 
Let $\pi:\CX\to B$ be an integral model of $X$ over $B$; namely, $\CX$ is a projective variety over $k$, and $\pi:\CX\to B$ is $k$-morphism whose generic fiber is isomorphic to $X\to \Spec K$. 
Let $\CL$ be a $\pi$-nef line bundle on $\CX$, i.e. $\CL$ has a non-negative degree on every projective curve in the fibers of $\pi$.  
Let $\CM$ be an ample line bundle on $B$. 
This defines a height function 
\[
h^\CM_\CL:X(\ol K)\longrightarrow\RR,\quad x\longmapsto \frac{1}{\deg(x)}\bar x\cdot \CL\cdot(\pi^*\CM)^{d-1}.
\]
Here $\bar x$ is the Zariski closure of $x$ in $\CX$. 
If $\CL$ is nef on $\CX$, then $h^\CM_\CL(x)\geq 0$ for every $x\in C(\ol K)$. 

The \emph{essential minimum} of $h_{\CL}^\CM$ is defined as 
\[
e_1(h_{\CL}^\CM)
=\sup_{U\subset X}\inf_{x\in U(\overline K)}h_{\CL}^\CM(x).
\]
Here the supremum is taken over all open subvarieties $U$ of $X$. 
The \emph{absolute minimum} of $h_{\CL}^\CM$ is defined as
\[
e_2(h_{\CL}^\CM)=\inf_{x\in X(\ol K)}h_{\CL}^\CM(x).
\]

Assume that the generic fiber $\CL_K=\CL|_X$ has a positive degree. Then the height of $X$ is defined by
$$
h^\CM_\CL(X)=\frac{1}{2\deg(\CL_K)} \CL^2\cdot(\pi^*\CM)^{d-1}.
$$

If $d=1$, then $K$ is the function field of one variable, and $\CM$ does not appear in the above definitions, so we drop the dependence on $\CM$ in $h^\CM_\CL$ and $e_i(h_{\CL}^\CM)$. The situation is parallel to the case of number fields treated by Zhang \cite{Zha92}. Then a geometric analogue of Zhang's theorem of successive minima in \cite[Thm. 6.3]{Zha92} asserts that
$$
e_1(h_{\CL}) \geq h_\CL(X) \geq \frac{1}{2}(e_1(h_{\CL})+e_2(h_{\CL})). 
$$
This holds if $d=1$, $\CL$ is $\pi$-nef, and $\deg(\CL_K)>0$, and the proof is similar to the loc. cit.. 

If $d>1$, it is natural to speculate
$$
e_1(h_{\CL}^\CM) \geq h^\CM_\CL(X) \geq \frac{1}{2}(e_1(h_{\CL}^\CM)+e_2(h_{\CL}^\CM)),
$$
assuming $\CL$ is $\pi$-nef, $\deg(\CL_K)>0$, and $\CM$ is ample on $B$.
In fact, the first inequality is claimed in \cite[Lem. 4.1]{Gub07}, and the second inequality is claimed in \cite[Prop. 4.3]{Gub07}. 
In the following, we prove that the second inequality is wrong by a counterexample.  

To introduce the counterexample, 
return to the situation of $C$ and $J$. 
Namely, let $C$ be a smooth projective curve  of genus $g>1$ over an algebraically closed field $k$, and let $J$ be the Jacobian variety of $C$. 

To utilize the above terminology of successive minima, we take $\pi:\CX\to B$ to be $p_2:C\times J\to J$.
Then $X=C_K$ for $K=k(J).$
Take 
$$\CM=\theta \in \Pic(J)$$ 
and 
\[
\CL=g\cdot \alpha_1+ \theta_2+Q\in\Pic(C\times J). 
\]
Here for simplicity, we write $\theta$ for the line bundle $\CO(\theta)$ in the Picard group. 
Then we have the following computational results, which disprove the second inequality speculated above. 

\begin{theorem}[Theorem \ref{second main}]
\label{second main thm}
Assume that $\rho(C\times J)=3$. Then 
\[
e_1(h_{\CL}^\theta)=e_2(h_{\CL}^\theta)=(g-\frac{1}{g})\cdot (g-1)!,
\]
and
\[
h_{\CL}^\theta(C_K)=(g-1)\cdot (g-1)!.
\]
Moreover, there is an infinite sequence $(x_n)_n$ in $C_K(\overline K)$ with $\deg(x_n)\to \infty$ and with 
$$
h^\theta_\CL(x_n) = (g-\frac{1}{g})\cdot (g-1)!.
$$
\end{theorem}

To prove Theorem \ref{second main thm}, the key is that in the height formula 
\[
h^\theta_\CL(x)= \frac{1}{\deg(x)}\bar x\cdot \CL\cdot(\pi^*\theta)^{d-1},
\]
the term $\bar x$ is an effective divisor on $C\times J$. 
Then we can use Theorem \ref{first main thm} to express $\bar x$ in terms $\alpha_1, \theta_2, Q$, and estimate the height by explicit calculations.

In the end, we make a few remarks on the truth of the inequalities
$$
e_1(h_{\CL}^\CM) \geq h^\CM_\CL(X) \geq \frac{1}{2}(e_1(h_{\CL}^\CM)+e_2(h_{\CL}^\CM)). 
$$
Here we assume that  $\CL$ is $\pi$-nef, $\deg(\CL_K)>0$, and $\CM$ is ample on $B$.
As mentioned above, both inequalities hold for $\dim B=1$. 
The first inequality always holds for general $\dim B$, but our counterexample disproves the second inequality for $\dim B=2$.

Let us examine the idea of ``generic curve'' in \cite[Prop. 5.11]{Gub08}. 
The process does not change the height of an algebraic point, but it brings more algebraic points. It does not change $h^\CM_\CL(X)$, but might decrease $e_1(h_{\CL}^\CM)$ and $e_2(h_{\CL}^\CM)$. 
Consequently, it reduces the first inequality from $\dim B>1$ to  $\dim B=1$, but it does not work for the second inequality. 
Alternatively, we can also apply the argument of \cite[Lem. 3.7]{XY22} to prove the first inequality, where we can choose a different type of ``generic curve'' such that $h^\CM_\CL(X)$ does not change.
This process does not bring more algebraic points, but it might decrease the height of an algebraic point, and thus might still decrease $e_1(h_{\CL}^\CM)$ and $e_2(h_{\CL}^\CM)$. 

Now we consider the arithmetic setting, where the projective varieties over $k$ are replaced by projective arithmetic varieties over $\ZZ$ and the line bundles are replaced by hermitian line bundles. 
If $\dim B=1$, this is Zhang's original theorem in \cite{Zha92}. 
If $\dim B>1$, the height $h_{\CL}^\CM$ is the Moriwaki height originally introduced by Moriwaki \cite{Mor00}, and we refer to \cite[Cor. 5.2]{Mor00} for a weaker version of the inequalities. 
With some effort, we could prove the first inequality for all $\dim B$ and disprove the second inequality for $\dim B>1$. 

Finally, we refer to Zhang \cite[Thm. 1.10]{Zha95} and Yuan--Zhang \cite[\S5.3]{YZ26} for adelic versions of the inequalities.

\subsubsection*{Acknowledgments}
The authors are grateful to Walter Gubler and Junyi Xie for communications related to this paper. 
The authors are also indebted to the anonymous referee for providing a reference for Lemma 3.1.

The authors would like to thank the support of the China--Russia Mathematics Center. The second author is supported by grants NO. 12250004 and NO. 12321001
from the National Science Foundation of China, and by the Xplorer Prize from the New Cornerstone Science Foundation.

\section{The cones}

Throughout this paper, let $k$ be an algebraically closed field. Let $C$ be a smooth projective curve over $k$ of genus $g>1$. Let $J$ be the Jacobian variety of $C$. 

\subsection{Preliminary results}

Denote by $$s:J\times J\lra J$$ the addition law of the abelian variety. 
For each $x \in J(k)$, denote by 
$$
T_x : J \longrightarrow J, \quad y \longmapsto y + x,
$$
the translation morphism by $x$. 

Let $\alpha$ be a line bundle on $C$ such that $(2g-2)\alpha\cong\omega_{C/k}$.
Consider the Abel-Jacobi map
\[
i_\alpha:C\longrightarrow J,\quad x\longmapsto x-\alpha.
\]
Let $\theta$ be the image of the morphism
\[
C^{g-1}\longrightarrow J,\quad (x_1,\cdots,x_{g-1})\longmapsto x_1+\cdots+x_{g-1}-(g-1)\alpha.
\]
It is well known that, by this special choice of $\alpha$, the divisor $\theta$ is ample on $J$ with a symmetric linear equivalence class. For the symmetry, we refer to see \cite[p. 74, eq. (1)]{Ser89}. 
We have a principal polarization 
\[
\phi:J\stackrel{\thicksim}\longrightarrow J^\vee,\quad  x\longmapsto 
T^*_{x}\theta-\theta.
\]
Note that our definition differs from the definition of $\phi_\theta$ in the first paragraph of \cite[p. 76]{Ser89} by a sign, since the loc. cit. takes $T_{x,*}\theta-\theta$ instead of $T^*_{x}\theta-\theta$.

Let $P$ be the line bundle on $J\times J$ which is the pull-back of the usual Poincar\'e bundle on $J\times J^\vee$ via 
\[
\id\times \phi:J\times J\longrightarrow J\times J^\vee.
\]
Then 
\[
P=s^*\theta-q_1^*\theta-q_2^*\theta,
\]
where $q_i:J\times J\to J$ is the projection to i-th factor.
See \cite[p. 76, (3)]{Ser89}.
Denote 
\[
Q=(i_\alpha\times \id)^*P,
\]
which is a line bundle on $C\times J$.
It is well-known that $Q$ has the following universal property for the Picard functor $\underline{\Pic}_{C/k}^0$. 

\begin{lemma}\label{computing Q}
For $x\in C(k)$ and $y\in J(k)$, we have line bundles
\[
Q|_{x\times J}= \phi(x-\alpha),\quad Q|_{C\times y}= y.    
\]
Here on the right-hand side of the last equality, $y$ denotes the line bundle on $C$ represented by $y$. 
\end{lemma}
\begin{proof}
We have
\[
Q|_{x\times J}=P|_{\{x-\alpha\}\times J}=T^*_{x-\alpha}\theta-\theta=\phi(x-\alpha),
\]
and
\[
Q|_{C\times y}=i_\alpha^*(P|_{y\times J})=y.
\] 
\end{proof}

%
%
%

As before, denote
$$
\alpha_1=p_1^*\alpha,\quad \theta_2= p_2^*\theta.
$$
For $m,n\in \ZZ$, consider the morphism 
$$
f_{m,n}: C\times J\longrightarrow J,\quad (x,y)\longmapsto m(x-\alpha)+ny. 
$$
As $\theta$ is ample on $J$, the pull-back $f_{m,n}^*\theta$ is nef on $C\times J$. 
This gives lots of examples of nef line bundles. 
We have the following explicit expression. 

\begin{proposition} \label{pullback}
In $\Pic(C\times J)$, we have
\[
f_{m,n}^*\theta=gm^2\cdot \alpha_1+n^2\cdot \theta_2+mn\cdot Q.
\]
\end{proposition}

\begin{proof}
Denote 
$$
F_{m,n}: J\times J\longrightarrow J,\quad (x,y)\longmapsto mx+ny. 
$$
It is the composition of the component-wise multiplication 
$$[m,n]:J\times J\longrightarrow J\times J, \quad (x,y)\longmapsto(mx,ny)$$
with the addition $s:J\times J\to J$.
It follows that 
$$
F_{m,n}^*\theta=[m,n]^*(s^*\theta)
=[m,n]^*(q_1^*\theta+q_2^*\theta+P)
=m^2\cdot q_1^*\theta+n^2\cdot q_2^*\theta+mn \cdot P.
$$
As 
$f_{m,n}: C\times J\to J$
is the composition of 
$i_\alpha\times \id: C\times J\to J\times J$
with 
$F_{m,n}: J\times J\to J$, we  have 
$$
f_{m,n}^*\theta=(i_\alpha\times \id)^*(F_{m,n}^*\theta)
=(i_\alpha\times \id)^*(m^2\cdot q_1^*\theta+n^2\cdot q_2^*\theta+mn\cdot P)
=m^2\cdot p_1^*(i_\alpha^*\theta)+n^2\cdot \theta_2+mn\cdot Q.
$$
By \cite[p. 75, (2)]{Ser89}, we have $i_\alpha^*\theta=g\alpha$. 
This finishes the proof. 
\end{proof}

In the following, we introduce a result, which is well known to experts.  
For example, some relevant references are \cite[\S VI.4]{Lan83},  \cite[Lem. 2.2.1]{Zha10}, and \cite[\S A-2]{MNP13}.
However, we are not able to find an exact reference, so we include a proof here. 

\begin{theorem}\label{decomposition}
There is a canonical isomorphism 
$$
\Psi:\Pic(C)\oplus \Pic(J) \oplus \End(J) \lra \Pic(C\times J)
$$
induced by the pull-back 
$p_1^*:\Pic(C)\to \Pic(C\times J)$, the pull-back  
$p_2^*: \Pic(J)\to \Pic(C\times J)$,
 and the map  
$$
\Psi_0:\End(J) \lra \Pic(C\times J),
\quad 
f \longmapsto (\mathrm{id}, f)^* Q. 
$$
Moreover, $\Psi$ induces a canonical isomorphism 
$$
\overline\Psi:
\NS(C)\oplus \NS(J) \oplus \End(J)\lra 
\NS(C\times J).
$$
\end{theorem}
\begin{proof}
We first prove that $\Psi$ is an isomorphism.
By definition, $J$ represents the Picard functor $\underline{\Pic}_{C/k}^0$. 
This gives  identities
$$
\underline{\Pic}_{C/k}^0(J)=\mathrm{Mor}_k(J, J)= \End(J) \oplus J(k).
$$
Here $\mathrm{Mor}_k(J, J)$ denotes the set of all $k$-morphisms 
$f:J\to J$,  $\End(J)$ is the subset of all $k$-morphisms 
$f:J\to J$ with $f(0)=0$, and $J(k)$ is the subset of all $k$-morphisms 
$f:J\to J$ with $f(J)=y$ for some $y\in J(k)$. 
By \cite[\S 8.1, Prop. 4]{BLR90}, we  have 
$$
\underline{\Pic}_{C/k}(J)=\Pic(C\times J)/p_2^*\Pic(J). $$
It follows that
$$
\underline{\Pic}_{C/k}^0(J)=\Pic(C\times J)_0/p_2^*\Pic(J),
$$
where
$\Pic(C\times J)_0$
is the subgroup of elements of $\Pic(C\times J)$ with degree zero on every fiber of $p_2:C\times J\to J$. 
Thus we have a canonical isomorphism 
$$
\Psi_1:\mathrm{Mor}_k(J, J) \lra \Pic(C\times J)_0/p_2^*\Pic(J).
$$
As $Q$ is the universal line bundle for the functor 
$\underline{\Pic}_{C/k}^0$, the map $\Psi_1$ is actually given by 
$$
\Psi_1: f\longmapsto (\mathrm{id},f)^* Q. 
$$

By Lemma \ref{computing Q}, 
for any $y\in J(k)$, we have 
$$
\Psi_1(y)= (\mathrm{id},y)^* Q=p_1^*(Q|_{C\times y})=p_1^*(y). 
$$
Then we have 
$$\Psi_1(J(k))=p_1^*\Pic^0(C).$$ 
So
$\Psi_1$ induces an isomorphism 
$$
\Psi_2:\End(J) \lra \Pic(C\times J)_0/(p_1^*\Pic^0(C)\oplus p_2^*\Pic(J))
$$
Composing with the natural isomorphism 
$$
\Pic(C\times J)_0/(p_1^*\Pic^0(C)\oplus p_2^*\Pic(J))\lra \Pic(C\times J)/(p_1^*\Pic(C)\oplus p_2^*\Pic(J)),
$$
we obtain an isomorphism 
$$
\Psi_3:\End(J) \lra \Pic(C\times J)/(p_1^*\Pic(C)\oplus p_2^*\Pic(J)). 
$$

Note that $\Psi_3$ is compatible with 
$$
\Psi_0:\End(J) \lra \Pic(C\times J). 
$$
Then the isomorphism $\Psi_3$ implies that $\Psi$ is also an isomorphism. 

Now we prove that $\overline\Psi$ is an isomorphism. 
The isomorphism
$\Psi$
implies an isomorphism 
$$
\underline{\Pic}_{C\times J/k}
=\underline{\Pic}_{C/k} \times \underline{\Pic}_{J/k} \times \End(J)
$$
of Picard schemes. 
Here $\End(J)$ is viewed as a constant group scheme over $k$. 
Then the identity component of $\underline{\Pic}_{C\times J/k}$ is 
$$
\underline{\Pic}_{C\times J/k}^0
=\underline{\Pic}_{C/k}^0 \times \underline{\Pic}_{J/k}^0.
$$
This gives a canonical isomorphism
$$
\Pic^0(C\times J)
=\Pic^0(C) \times \Pic^0(J).
$$
By taking quotients, the isomorphism  
$\Psi$ induces the isomorphism $\overline \Psi$. 
\end{proof}

\subsection{Determination of the cones}

With the above preparation, we are ready to have Theorem \ref{first main thm}. 
For convenience, we duplicate the theorem as follows.

\begin{theorem}[Theorem \ref{first main thm}]\label{cone thm}
Assume that $\rho(C\times J)=3$. Then
	\[
	\Amp(C\times J)=\BBig(C\times J)=\{a\cdot \alpha_1+b\cdot \theta_2+c\cdot Q:a>0,\ b> 0, \ ab> gc^2\}
	\]
and
	\[
	\Nef(C\times J)=\Psef(C\times J)=\{a\cdot \alpha_1+b\cdot \theta_2+c\cdot Q:a\geq0,\ b\geq 0, \ ab\geq gc^2\}.
	\]
\end{theorem}

\begin{proof}
Denote 
\[
\Sigma= \{a\cdot \alpha_1+b\cdot \theta_2+c\cdot Q:a\geq0,\ b\geq 0, \ ab\geq gc^2\}.
\]
It suffices to prove 
\[
\Psef(C\times J) \ \subseteq \
\Sigma
\ \subseteq \
\Nef(C\times J). 
\]
In fact, this forces equalities by $\Nef(C\times J) \ \subseteq 
\Psef(C\times J).$
The interiors give the expression for $\Amp(C\times J)$ and $\BBig(C\times J)$. 

We first prove $\Sigma \subseteq  \Nef(C\times J)$. 
By Proposition \ref{pullback}, the line bundle
\[
L_{m,n}=gm^2\cdot \alpha_1+n^2\cdot \theta_2+mn\cdot Q
\]
is nef on $C\times J$ for any $m,n\in\ZZ$. 
By homogeneity, $L_{m,n}$ is nef for all  $m,n\in\QQ$.
By approximation, $L_{m,n}$ is nef for all $m,n\in\RR$.
Moreover, $L_{m,n}+t_1\alpha_1+t_2\theta_2$ is nef for all $m,n\in\RR$ and all $t_1,t_2\geq 0$. 
This proves that 
$\Sigma \subseteq  \Nef(C\times J).$

Now we prove $\Psef(C\times J) \subseteq \Sigma$. 
By taking the closure, it suffices to prove $\BBig(C\times J) \subseteq \Sigma$. 
Let 
$$L=a\cdot \alpha_1+b\cdot \theta_2+c\cdot Q$$ 
be an element of $\BBig(C\times J)$ with $a,b,c\in \RR$. 
We first claim that $a,b>0$. 
In fact, by Lemma \ref{computing Q},
\[
L|_{x\times J}=b\cdot \theta+c \cdot \phi(x-\alpha), \quad x\in C(k)
\]
and 
\[
L|_{C\times y}=a\cdot \alpha+c \cdot y, \quad y\in J(k).
\]
By the Kodaira lemma (cf. \cite[Cor. 2.2.7]{Laz04}), we can write $L=A+E$ for an ample class $A$ and an effective class $E$ on $J$. 
It follows that $L|_{x\times J}$ (resp. $L|_{C\times y}$) is big as long as $x\times J$ 
(resp. ${C\times y}$) is not contained in the support of $E$. 
It follows that $\deg(L|_{C\times y})>0$, and thus $a>0$. 
On the other hand, the bigness gives $L|_{x\times J}\cdot \theta^{g-1}>0$, which implies $b>0$ as $\phi(x-\alpha)$ is algebraically equivalent to 0.
These prove $a,b>0$. 

It remains to prove $ab\geq gc^2$. 
Write 
$$
L=L'+ (a-gc^2/b)\alpha_1 
$$
with 
$$L'=(gc^2/b)\cdot \alpha_1+b\cdot \theta_2+c\cdot Q.$$
Note that $L'=L_{m,n}$ for some $m,n\in\RR_{>0}$. 
We claim that 
$$
L'^{g+1}=0, \quad L'^g\cdot \alpha_1=b^g\cdot g!.
$$
The first equality holds if $m,n\in \ZZ$ by Proposition \ref{pullback} by dimension reasons, holds for $m,n\in \QQ$ by homogeneity, and thus holds for all $m,n\in \RR$ by continuity. 
For second equality, taking a point $x\in C(k)$, we have 
$$
L'^g\cdot \alpha_1
=L'^g\cdot (x\times J)
=(L'|_{x\times J})^g
=\big(b\cdot \theta+c\cdot \phi(x-\alpha)\big)^g
=\big(b \theta\big)^g
=b^g\cdot g!.
$$
Here the third equality follows from the first equality of Lemma \ref{computing Q}, and the fourth equality holds since  $\phi(x-\alpha)$ is algebraically equivalent to 0.

With these two equalities on intersection numbers, the nefness of $L'$ gives
$$
0\leq L\cdot L'^g =(a-gc^2/b)(b^g\cdot g!). 
$$
This finishes the proof. 
\end{proof}

Recall that a line bundle is semi-ample if some positive multiple is globally generated. 
As a surprising property from the proof, every element of $\Nef(C\times J)$ can be represented by a semi-ample element of $\Nef(C\times J)$. This eventually boils down to the fact that the boundary element  
\[
L_{m,n}=gm^2\cdot \alpha_1+n^2\cdot \theta_2+mn\cdot Q 
\]
is semi-ample. In fact, it is the pull-back of an ample class $\theta$.
As we will see, the last statement of Theorem \ref{second main thm} also relies on this property in some way.

\section{The successive minima}

The goal of this section is to prove Theorem \ref{second main thm}. 
Resume the notations in the previous section. 
On $C\times J$, we have the nef line bundle
\[
\CL=f_{1,1}^*\theta=g \alpha_1+\theta_2+Q.
\]
Consider the family 
\[
p_2:C\times J\longrightarrow J.
\]
Denote by $K=k(J)$ the function field of $J$. Recall that we have a  height function
\[
h_\CL^\theta: C_K(\ol K)\longrightarrow \RR,\quad x\longmapsto \frac{1}{\deg(x)}\bar x\cdot\CL\cdot \theta_2^{g-1}.
\]
We also have a height 
\[
h_\CL^\theta(C_K)=\frac{\CL^2\cdot \theta_2^{g-1}}{2 \deg(\CL_K)}.
\]
To compute them, we first need the following numerical result. 

\begin{lemma} \label{computing intersection}
For $a,b,c\in\RR$, 
\[
\left(a\cdot \alpha_1+b\cdot \theta_2+c\cdot Q\right)\cdot \CL\cdot \theta_2^{g-1}=(a+gb-2c)\cdot g!.
\]
\end{lemma}
\begin{proof}
By Lemma \ref{computing Q},  
$$
 \alpha_1\cdot Q \cdot \theta_2^{g-1}
=(Q|_{x\times J})\cdot (\theta_2|_{x\times J})^{g-1}
=0, 
$$
and 
$$
 \theta_2 \cdot Q \cdot \theta_2^{g-1}
=\deg(Q|_{C\times y}) (\theta^{g})
=0. 
$$
Here $x\in C(k)$ and $y\in J(k)$ are points. 
It follows that 
$$
\alpha_1\cdot \CL\cdot \theta_2^{g-1}
=\alpha_1\cdot \theta_2 \cdot \theta_2^{g-1}
=  \theta^{g}= g!, 
$$
and 
$$
\theta_2\cdot \CL\cdot \theta_2^{g-1}
= \CL\cdot \theta_2^{g}
= \deg(\CL|_{C\times y}) (\theta^{g})
=g\cdot g!.
$$

It remains to compute
$$
Q\cdot \CL\cdot \theta_2^{g-1}
= Q\cdot (g \alpha_1+\theta_2+Q)\cdot \theta_2^{g-1}
= Q^2\cdot \theta_2^{g-1}.
$$
We need the classical numerical equivalence relation
$$
\theta^{g-1} \equiv (g-1)!\, [C]
$$
of $1$-cycles on $J$.
If $k=\CC$, this follows from Poincar\'e's formula (cf. \cite[\S I.5]{ACGH85}) for cohomology classes;
in general, this follows from the recursive formula in \cite[\S2, (4)]{Mat62}.
It follows that 
$$
Q^2\cdot \theta_2^{g-1}
=(g-1)!\, Q^2 \cdot [C\times C]
=(g-1)! (j^*Q)^2. 
$$
Here $j: C\times C\to C\times J$ is the morphism sending $(x,y)$ to $(x,i_\alpha(y))$.

We claim that 
$$j^*Q=\Delta-q_1^* \alpha- q_2^* \alpha$$
in $\Pic(C\times C)$, where $q_1, q_2:C\times C\to C$ denote the projections, and $\Delta$ denotes the diagonal of $C\times C$.
To prove the claim, by the see-saw theorem (cf. \cite[p. 54, Cor. 6]{Mum08}), it suffices to prove that for any $x,y\in C(k)$, 
$$
(j^*Q)|_{x\times C}=(\Delta-q_1^* \alpha- q_2^* \alpha)|_{x\times C},\quad
(j^*Q)|_{C\times y}=(\Delta-q_1^* \alpha- q_2^* \alpha)|_{C\times y}. 
$$
Both results follow from Lemma \ref{computing Q}. 

With the claim, we easily have 
$$
(j^*Q)^2=(\Delta-q_1^* \alpha- q_2^* \alpha)^2
=(\Delta-q_1^* \alpha- q_2^* \alpha)\cdot \Delta
=-2g.
$$ 
Here $\Delta^2=2-2g$ by $\CO(-\Delta)|_{\Delta}= \omega_{\Delta/k}$. 
It follows that
$$Q^2\cdot \theta_2^{g-1}=-2\cdot (g-1)!.$$
The proof is complete.
\end{proof}

Now we are ready to prove Theorem \ref{second main thm}. 

\begin{theorem}[Theorem \ref{second main thm}]
\label{second main}
Assume that $\rho(C\times J)=3$. Then 
\[
e_1(h_{\CL}^\theta)=e_2(h_{\CL}^\theta)=(g-\frac{1}{g})\cdot (g-1)!,
\]
and
\[
h_{\CL}^\theta(C_K)=(g-1)\cdot(g-1)!.
\]
Moreover, there is an infinite sequence $(x_n)_n$ in $C_K(\overline K)$ with $\deg(x_n)\to \infty$ and with 
$$
h^\theta_\CL(x_n) = (g-\frac{1}{g})\cdot (g-1)!.
$$
\end{theorem}

\begin{proof}
By Lemma \ref{computing intersection}, we easily have
\[
h_{\CL}^\theta(C_K)=\frac{\CL^2\cdot \theta_2^{g-1}}{2\deg(\CL_K)}=\frac{(2g-2)\cdot g!}{2g}=(g-1)\cdot(g-1)!.
\]
Now we compute the successive minima.

For any point $x\in C_K(\ol K)$, its Zariski closure $\bar x$ is an effective divisor on $C\times J$. Suppose that in $\NS(C\times J)_\RR$, we have
\[
\bar x\equiv a\cdot \alpha_1+b\cdot\theta_2+c\cdot Q, \quad a,b,c\in \QQ.
\]
We have 
$$
\deg(x)= \deg(\bar x|_{C\times y})=a>0. 
$$
By Theorem \ref{cone thm} for the pseudo-effective cone, we have $b\geq 0$, and $ab\geq gc^2$. 
Write 
$s=b/a$ and $t=c/a$. Then
$$
\bar x\equiv a( \alpha_1+s\cdot\theta_2+t\cdot Q), \quad s\geq  gt^2.
$$
By Lemma \ref{computing intersection}, we have 
\begin{align*}
h^\theta_\CL(x)
&=\frac{1}{\deg(x)}\bar x\cdot\CL\cdot \theta_2^{g-1}\\
&=( \alpha_1+s\cdot\theta_2+t\cdot Q)\cdot\CL\cdot \theta_2^{g-1}\\
&= (1+gs-2t)\cdot g!  \\
&\geq (1+g^2t^2-2t)\cdot g!\\
&\geq (1-\frac{1}{g^2})\cdot g!\\
&= (g-\frac{1}{g})\cdot (g-1)!.
\end{align*}
As a consequence, 
\[
e_1(h_{\CL}^\theta)\geq e_2(h_{\CL}^\theta)\geq (g-\frac{1}{g})\cdot (g-1)!. 
\]
Note that the last statement of the theorem implies 
\[
e_1(h_{\CL}^\theta)\leq (g-\frac{1}{g})\cdot (g-1)!,
\]
which forces equalities.

It remains to prove the last statement of the theorem. 
Consider the condition of equality in the above proof of 
$$
h^\theta_\CL(x)
\geq (g-\frac{1}{g})\cdot (g-1)!, 
$$ 
we have to take 
$$t= \frac{1}{g^2},\quad
s= \frac{1}{g^3}.$$
Hence, we take 
$$
D=\alpha_1+\frac{1}{g^3}\cdot\theta_2+\frac{1}{g^2}\cdot Q.
$$
Note that this divisor class lies in the boundary of the nef cone of $C\times J$. 
By Proposition \ref{pullback}, we have  
$$
f_{g,1}^* \theta = g^3\cdot \alpha_1+ \theta_2+g\cdot Q
=g^3 D.
$$

As $\theta$ is ample on $J$, the multiple $n\theta$ is very ample for sufficiently large $n$. 
By Bertini's theorem, such $n\theta$ is linearly equivalent to a smooth (and connected) divisor $H_n$ on $J$. We will consider its pull-back via 
$$
f_{g,1}: C\times J\longrightarrow J,\quad (x,y)\longmapsto g(x-\alpha)+y. 
$$ 
The fiber of $f_{g,1}$ above every $z\in J(k)$ is isomorphic to $C$, and thus is 1-dimensional and smooth. By the miracle flatness theorem (cf. \cite[Thm. 23.1]{Mat}),  $f_{g,1}$ is flat, and thus it is further smooth. 

As a consequence, the inverse image $E_n= f_{g,1}^{-1}(H_n)$ is smooth over $H_n$.
In particular, $E_n$ is also smooth and connected.   
In $\Pic(C\times J)$, we have 
$$
E_n=f_{g,1}^*(H_n)=f_{g,1}^*(n\theta)=g^3 n D. 
$$

The divisor $E_n$ restricts to a closed point on the generic fiber $C_K$ and thus induces an algebraic point $x_n\in C_K(\overline K)$. 
The degree 
$$
\deg(x_n)= \deg(E_n |_{C\times y})=g^3 n \deg(D |_{C\times y})=g^3 n. 
$$
We also have 
\[
h^\theta_\CL(x_n)
=\frac{1}{\deg(x_n)} E_n\cdot \CL\cdot \theta_2^{g-1}
=\frac{1}{g^3n} g^3nD\cdot \CL\cdot \theta_2^{g-1}
=D\cdot \CL\cdot \theta_2^{g-1}
=(g-\frac{1}{g})\cdot (g-1)!.
\]
Here the last equality follows from Lemma \ref{computing intersection} again. 
This finishes the proof. 
\end{proof}

\

{\footnotesize
\noindent Ruoyi Guo

\noindent Address: \emph{School of Mathematical Sciences, Peking University, Haidian District, Beijing 100871, China}

\noindent  Email: \emph{guoruoyi@pku.edu.cn}

\

\noindent Xinyi Yuan

\noindent Address: \emph{BICMR, Peking University, Haidian District, Beijing 100871, China}

\noindent Email: \emph{yxy@bicmr.pku.edu.cn}
}

 \end{document}